\documentclass{amsart}
\usepackage{amsmath, amsthm, amsfonts, mathtools, tikz-cd, amssymb}
\usepackage{enumitem}

\makeatletter
\def\paragraph{\@startsection{paragraph}{4}%
  \z@\z@{-\fontdimen2\font}%
  {\normalfont\bfseries}}
\makeatother

\title[Topological complexity of toral relatively hyperbolic groups]{On the topological complexity of toral relatively hyperbolic groups}
\author{Kevin Li}
\address{School of Mathematical Sciences, University of Southampton, Southampton SO17 1BJ, United Kingdom}
\email{kevin.li@soton.ac.uk}
\date{\today}

\theoremstyle{definition}
\newtheorem{defn}{Definition}[section]
\newtheorem{setup}[defn]{Setup}
\theoremstyle{plain}
\newtheorem{thm}[defn]{Theorem}
\newtheorem{lem}[defn]{Lemma}

\newtheorem{cor}[defn]{Corollary}
\theoremstyle{remark}

\newcommand{\IZ}{\ensuremath\mathbb{Z}}
\newcommand{\F}{\ensuremath\mathcal{F}}
\newcommand{\E}{\ensuremath\mathcal{E}}
\newcommand{\TR}{\ensuremath\mathcal{TR}}
\newcommand{\FIN}{\ensuremath\mathcal{FIN}}
\newcommand{\VCY}{\ensuremath\mathcal{VCY}}
\newcommand{\ALL}{\ensuremath\mathcal{ALL}}
\newcommand{\calH}{\ensuremath\mathcal{H}}
\newcommand{\calP}{\ensuremath\mathcal{P}}
\newcommand{\calM}{\ensuremath\mathcal{M}}
\newcommand{\calD}{\ensuremath\mathcal{D}}
\newcommand{\OFG}{{\ensuremath\mathcal{O}_\mathcal{F}G}}
\newcommand{\EFG}{\ensuremath E_\mathcal{F}G}
\newcommand{\EEG}{\ensuremath E_\mathcal{E}G}
\newcommand{\Mod}{\ensuremath\text{-}\mathbf{Mod}}
\newcommand{\enum}{\rm{(\roman*)}}
\newcommand{\spann}[1]{{\ensuremath \langle{#1}\rangle}}
\DeclareMathOperator{\Hom}{Hom}
\DeclareMathOperator{\id}{id}
\DeclareMathOperator{\cd}{cd}
\DeclareMathOperator{\TC}{TC}
\DeclareMathOperator{\hdim}{hdim}

\begin{document}
\maketitle

\begin{abstract}
	We prove that the topological complexity $\TC(\pi)$ equals $\cd(\pi\times\pi)$ for certain toral relatively hyperbolic groups $\pi$.
\end{abstract}

\section{Introduction}

The (reduced) topological complexity $\TC(X)$ of a space $X$ is defined as the minimal integer $n$ for which there exists a cover of $X\times X$ by $n+1$ open subsets $U_0,\ldots,U_n$ such that the path fibration $X^{[0,1]}\to X\times X$ admits a local section over each $U_i$. This quantity, which is similar in spirit to the classical Lusternik--Schnirelmann category, was introduced by Farber~\cite{Farber03} in the context of robot motion planning. In fact, $\TC(-)$ is a homotopy invariant and hence one can define the topological complexity $\TC(\pi)$ of a group $\pi$ to be $\TC(B\pi)$, where $B\pi$ is the classifying space for $\pi$. There are bounds $\cd(\pi)\le \TC(\pi)\le \cd(\pi\times\pi)$, where $\cd(-)$ denotes the cohomological dimension. However, the precise value of $\TC(\pi)$ is known only for a small class of groups, which contains for instance the abelian groups, hyperbolic groups, free products of the form $H\ast H$ for $H$ geometrically finite, right-angled Artin groups, and certain subgroups of braid groups. We refer to~\cite{Farber-Mescher20} and~\cite{Dranishnikov20} for a more thorough account on this topic.

It is the decisive insight of~\cite{Farber19} that the topological complexity of groups can be expressed in terms of classifying spaces for families of subgroups, which are well-studied objects in equivariant topology. For a family $\F$ of subgroups of a group $G$, the classifying space $\EFG$ is a terminal object, up to $G$-homotopy, among $G$-CW-complexes with stabilizers in $\F$.
Farber, Grant, Lupton, and Oprea showed that $\TC(\pi)$ equals the minimal integer $n$ for which the canonical $(\pi\times\pi)$-map $E(\pi\times\pi)\to E_\calD(\pi\times\pi)$ is equivariantly homotopic to a map with values in the $n$-skeleton $E_\calD(\pi\times\pi)^{(n)}$. Here $\calD$ is the family of subgroups of $\pi\times\pi$ consisting of all conjugates of the diagonal subgroup $\Delta(\pi)$ and their subgroups. Using this characterization of $\TC(\pi)$, in a recent breakthrough Dranishnikov~\cite{Dranishnikov20} has computed the topological complexity of torsionfree hyperbolic groups and more generally, of geometrically finite groups with cyclic centralizers.

\begin{thm}[Dranishnikov]\label{thm:Dranishnikov}
	Let $\pi$ be a geometrically finite group with $\cd(\pi)\ge 2$ such that the centralizer $Z_\pi(b)$ is cyclic for any $b\in\pi\setminus\{e\}$. Then $\TC(\pi)=\cd(\pi\times\pi)$.
\end{thm}
Recall that a group $\pi$ is called \emph{geometrically finite} if it admits a finite model for $B\pi$.
Note that for geometrically finite groups $\pi$ we have $\cd(\pi\times\pi)=2\cd(\pi)$, see~\cite{Dranishnikov19}.
Previously, Farber and Mescher~\cite{Farber-Mescher20} had shown for groups $\pi$ as in Theorem~\ref{thm:Dranishnikov} that $\TC(\pi)$ equals either $\cd(\pi\times\pi)$ or $\cd(\pi\times\pi)-1$.
The main contribution of the present note is the following generalization of Theorem~\ref{thm:Dranishnikov}.

\begin{thm}\label{thm:main thm}
	Let $\pi$ be a torsionfree group with $\cd(\pi)\ge 2$. Suppose that $\pi$ admits a malnormal collection of abelian subgroups $\calP=\{P_i\ |\ i\in I\}$ satisfying $\cd(P_i\times P_i)<\cd(\pi\times\pi)$ such that the centralizer $Z_\pi(b)$ is cyclic for any $b\in\pi$ that is not conjugate into any of the $P_i$. 
	Then $\TC(\pi)=\cd(\pi\times\pi)$.
\end{thm}

Recall that a set $\calP=\{P_i\ |\ i\in I\}$ of subgroups of $\pi$ is called a \emph{malnormal collection} if for any $P_i,P_j\in\calP$ and $g\in\pi$, we have $gP_ig^{-1}\cap P_j=\{e\}$ or $i=j$ and $g\in P_i$.
Our main examples of groups satisfying the assumptions of Theorem~\ref{thm:main thm} are torsionfree relatively hyperbolic groups $\pi$ with $\cd(\pi)\ge 2$ and finitely generated abelian peripheral subgroups $P_1,\ldots,P_k$ satisfying $\cd(P_i)<\cd(\pi)$.
Note that Theorem~\ref{thm:main thm} recovers Theorem~\ref{thm:Dranishnikov} as a special case when $\calP$ consists only of the trivial subgroup and that the assumption of geometric finiteness has been dropped.

In light of the upper bound $\TC(\pi)\le \cd(\pi\times\pi)$, Theorem~\ref{thm:Dranishnikov} and Theorem~\ref{thm:main thm} are statements about the maximality of topological complexity. They share a common strategy of proof based on the characterization of $\TC(\pi)$ in terms of classifying spaces from~\cite{Farber19}. Namely, we construct a ``small" model for $E_\calD(\pi\times\pi)$ from $E(\pi\times\pi)$ allowing us to show that the map $E(\pi\times\pi)\to E_\calD(\pi\times\pi)$ induces a non-trivial map on cohomology in degree $\cd(\pi\times\pi)$. Hence one has equality $\TC(\pi)=\cd(\pi\times\pi)$.
Nevertheless, even for the case when $\calP$ consists only of the trivial subgroup, our proof is different from Dranishnikov's. He constructed a specific model for $E_\calD(\pi\times\pi)$ and used cohomology with compact support, while we employ a general construction due to L\"uck and Weiermann and use equivariant Bredon cohomology. L\"uck and Weiermann's construction (Theorem~\ref{thm:LW}) is a general recipe to efficiently construct $\EFG$ from $\EEG$ for two families of subgroups $\E\subset \F$ of a group $G$ satisfying a certain maximality condition. While for the group $\pi\times\pi$ this condition is not satisfied for the families $\{\{e\}\}\subset\calD$, we define an intermediate family $\{\{e\}\}\subset\F_1\subset\calD$ such that we can apply two iterations of the construction.
\\
\paragraph{Acknowledgments.} 
	The present note is part of the author's PhD project under the supervision of Nansen Petrosyan and Ian Leary, who we thank for their support. We are grateful to Pietro Capovilla for interesting discussions about the paper~\cite{CLM20} and his master's thesis. We thank Mark Grant for helpful comments on an earlier version of this note.

\section{Preliminaries on classifying spaces for families}
We briefly review the notion of classifying spaces for families of subgroups due to tom Dieck and their equivariant Bredon cohomology. For a survey on classifying spaces for families we refer to~\cite{Lueck05survey} and for an introduction to Bredon cohomology to~\cite{Fluch10}. Let $G$ be a group, which shall always mean a discrete group. 

A \emph{family of subgroups $\F$} is a non-empty set of subgroups of $G$ that is closed under conjugation by elements of $G$ and under taking subgroups. Typical examples are $\TR=\{\{e\}\}$, $\FIN=\{\text{finite subgroups}\}$, $\VCY=\{\text{virtually cyclic subgroups}\}$, and $\ALL=\{\text{all subgroups}\}$. For a set $\calH$ of subgroups of $G$, one can consider $\F\spann{\calH}=\{\text{conjugates of subgroups in $\calH$ and their subgroups}\}$ which is the smallest family containing $\calH$ and called the \emph{family generated by $\calH$}. When $\calH=\{H\}$ consists of a single subgroup, we denote $\F\spann{\{H\}}$ instead by $\F\spann{H}$ and call it the \emph{family generated by $H$}. For a family $\F$ of subgroups of $G$ and any subgroup $H\subset G$, we denote by $\F|_H$ the family $\{K\cap H\ |\ K\in\F\}$ of subgroups of $H$. (In the literature this family is sometimes denoted by $\F\cap H$ instead.)

A \emph{classifying space $\EFG$ for the family $\F$} is a terminal object in the $G$-homotopy category of $G$-CW-complexes with stabilizers in $\F$. It can be shown that $\EFG$ always exists and that a $G$-CW-complex $X$ is a model for $\EFG$ if and only if the fixed-point set $X^H$ is contractible for $H\in\F$ and empty otherwise. In particular, there exists a $G$-map $EG\to \EFG$ which is unique up to $G$-homotopy.

The \emph{orbit category $\OFG$} has as objects $G/H$ for $H\in\F$ and as morphisms $G$-maps. 
Let $\OFG\Mod$ denote the category of contravariant functors $M\colon \OFG\to \IZ\Mod$ with values in the category of $\IZ$-modules, which are called \emph{$\OFG$-modules}.
For a $G$-CW-complex $X$ with stabilizers in $\F$, the \emph{$G$-equivariant Bredon cohomology $H^*_G(X;M)$} with coefficients in an $\OFG$-module $M$ is the cohomology of the cochain complex $\Hom_{\OFG\Mod}(C_*(X^?),M)$, where $C_*(X^?)(G/H)=C_*(X^H)$ is the cellular chain complex. 
\\
\paragraph{Passage to larger families}
Let $G$ be a group and $\E\subset\F$ be two families of subgroups.

We say that $G$ satisfies \emph{condition $(M_{\E\subset\F})$} if every element $H\in\F\setminus\E$ is contained in a unique element $M\in\F\setminus\E$ which is maximal in $\F\setminus\E$ (with respect to inclusion).
We say that $G$ satisfies \emph{condition $(NM_{\E\subset\F})$} if every maximal element $M\in\F\setminus\E$ is self-normalizing, i.e.\ $M$ equals its normalizer $N_GM$ in $G$.
Let $\calM=\{M_i\ |\ i\in I\}$ be a complete set of representatives for the conjugacy classes of maximal elements in $\F\setminus\E$, i.e.\ each $M_i$ is maximal in $\F\setminus \E$ and any maximal element in $\F\setminus\E$ is conjugate to precisely one of the $M_i$. The following~\cite[Corollary 2.8]{Lueck-Weiermann12} is a special case of a more general construction due to L\"uck and Weiermann.

\begin{thm}[L\"uck--Weiermann]\label{thm:LW}
	Let $G$ be a group satisfying condition $(M_{\E\subset \F})$ for two families of subgroups $\E\subset \F$.
	Consider a cellular $G$-pushout of the form
	\[\begin{tikzcd}
		\coprod_{i\in I} G\times_{N_GM_i}E_{\E|_{N_GM_i}}(N_GM_i)\ar{r}{\varphi} \ar{d}[swap]{\coprod_{i\in I}\id_G\times_{N_GM_i}f_i} & E_\E G\ar{d} \\
		\coprod_{i\in I} G\times_{N_GM_i}E_{\ALL|_{M_i}\cup \E|_{N_GM_i}}(N_GM_i)\ar{r} & X
	\end{tikzcd}\]
	such that each $f_i$ is a cellular $N_GM_i$-map and $\varphi$ is an inclusion of $G$-CW-complexes, or such that each $f_i$ is an inclusion of $N_GM_i$-CW-complexes and $\varphi$ is a cellular $G$-map.
	Then $X$ is a model for $\EFG$.
\end{thm}

Note that a $G$-pushout as in Theorem~\ref{thm:LW} with maps $f_i$ and $\varphi$ as required always exists by using equivariant cellular approximation and mapping cylinders.

\begin{cor}\label{cor:LW}
	Let $G$ be a group and $\E\subset\F$ be two families of subgroups.
\begin{enumerate}[label=\enum]
	\item\label{enum:TRIV} If $G$ satisfies condition $(M_{\TR\subset \F})$, then a model for $\EFG$ can be constructed as a $G$-pushout of the form
	\[\begin{tikzcd}
		\coprod_{i\in I} G\times_{N_GM_i}E(N_GM_i)\ar{r}\ar{d} & EG\ar{d} \\
		\coprod_{i\in I} G\times_{N_GM_i}E(N_GM_i/M_i) \ar{r} & \EFG \,;
	\end{tikzcd}\]

	\item\label{enum:NM} If $G$ satisfies conditions $(M_{\E\subset \F})$ and $(NM_{\E\subset \F})$, then a model for $\EFG$ can be constructed as a $G$-pushout of the form
	\[\begin{tikzcd}
		\coprod_{i\in I} G\times_{M_i}E_{\E|_{M_i}} M_i\ar{r}\ar{d} & E_\E G\ar{d} \\
		\coprod_{i\in I} G/M_i \ar{r} & \EFG \,.
	\end{tikzcd}\]
\end{enumerate}
	\begin{proof}
		This follows from Theorem~\ref{thm:LW} by observing that if $\E|_{N_GM_i}\subset \ALL|_{M_i}$, then a model for $E_{\ALL|_{M_i}\cup \E|_{N_GM_i}}(N_GM_i)$ is given by $E(N_GM_i/M_i)$ regarded as a $N_GM_i$-CW-complex.
	\end{proof}
\end{cor}

\paragraph{Homotopy dimension and cohomological dimension of maps}
Let $G$ be a group and $\E\subset\F$ be two families of subgroups. The following notation is not standard. 

We denote by $\hdim_{\E\subset\F}(G)$ the minimal integer $n$ for which the canonical $G$-map $E_\E G\to \EFG$ is $G$-homotopic to a $G$-map with values in the $n$-skeleton $(\EFG)^{(n)}$. We denote by $\cd_{\E\subset\F}(G)$ the maximal integer $k$ for which the induced map on Bredon cohomology $H^k_G(\EFG;M)\to H^k_G(\EEG;M)$ is non-trivial for some $\OFG$-module $M$. 
One clearly has the inequality
\[
	\cd_{\E\subset\F}(G)\le \hdim_{\E\subset\F}(G) \,.
\]

\paragraph{Topological complexity as homotopy dimension}
Let $\pi$ be a group and $\Delta(\pi)\subset\pi\times\pi$ be the diagonal subgroup. Consider the family $\calD\coloneqq \F\spann{\Delta(\pi)}$ of subgroups of $\pi\times\pi$ that is generated by $\Delta(\pi)$.
The following is the main result of~\cite[Theorem 3.3]{Farber19}.
\begin{thm}[Farber--Grant--Lupton--Oprea]\label{thm:FGLO}
	Let $\pi$ be a group. Then $\TC(\pi)=\hdim_{\TR\subset\calD}(\pi\times\pi)$.
\end{thm}
Theorem~\ref{thm:FGLO} was recently generalized to families generated by a single subgroup in~\cite[Theorem 1.1]{BCE20} and to arbitrary families in~\cite[Proposition 7.5]{CLM20}.

\section{Structure of the diagonal family of $\pi\times\pi$}\label{sec:diagonal family}

Let $\pi$ be a group and $\Delta\colon \pi\to \pi\times\pi$ be the diagonal map. 
For a subset $S\subset \pi$, denote by $Z_\pi(S)$ the centralizer of $S$ in $\pi$. 
The following notation is adopted from~\cite{Farber19} and~\cite{Dranishnikov20}.

For $\gamma\in\pi$ and a subset $S\subset\pi$, define the subgroup $H_{\gamma,S}$ of $\pi\times\pi$ to be
\[
	H_{\gamma,S} \coloneqq (\gamma,e)\cdot \Delta(Z_\pi(S))\cdot (\gamma^{-1},e) \,.
\]
When $S$ is a singleton set $\{b\}$, we write $H_{\gamma,b}$ instead of $H_{\gamma,\{b\}}$. Note that $H_{e,e}=\Delta(\pi)$. 
The proof of the following identities is elementary and left to the reader.

\begin{lem}\label{lem:HgS operations}
	Let $\gamma,\delta\in \pi$ and $S,T\subset \pi$ be subsets. Then the following hold:
	\begin{enumerate}[label=\enum]
		\item\label{enum:conjugation} $(g,h)\cdot H_{\gamma,S}\cdot (g^{-1},h^{-1})=H_{g\gamma h^{-1},hSh^{-1}}$ for any $(g,h)\in\pi\times\pi$;
		\item\label{enum:intersection} $H_{\gamma,S}\cap H_{\delta,T}=H_{\gamma,S\cup T\cup \{\delta^{-1}\gamma\}}$;
		\item\label{enum:normalizer} $N_{\pi\times\pi}H_{\gamma,S}=\{(\gamma kh\gamma^{-1},h)\in\pi\times\pi \ |\ h\in N_\pi(Z_\pi(S)), k\in Z_\pi(Z_\pi(S))\}$.
	\end{enumerate}
\end{lem}

We define the families $\F_1\subset\calD$ of subgroups of $\pi\times\pi$ to be
\begin{equation}\label{eqn:D F1}
\begin{aligned}
	\calD &\coloneqq \F\spann{\Delta(\pi)}\,;\\
	\F_1 &\coloneqq \F\spann{\{H_{\gamma,b}\ |\ \gamma\in\pi,b\in\pi\setminus\{e\}\}}\,.
\end{aligned}
\end{equation}
In view of Lemma~\ref{lem:HgS operations}~\ref{enum:conjugation} and~\ref{enum:intersection}, the family $\F_1$ is generated by the intersections of conjugates of the diagonal subgroup $\Delta(\pi)$.

\begin{lem}\label{lem:F1 D}
	Let $\pi$ be a group. Then condition $(M_{\F_1\subset\calD})$ holds for the group $\pi\times\pi$. Moreover, if the center $Z_\pi(\pi)$ of $\pi$ is trivial, then condition $(NM_{\F_1\subset\calD})$ holds.
	\begin{proof}
		If $\F_1$ equals $\calD$, then the statement is vacuous, so we may assume that $\F_1$ is strictly contained in $\calD$.
		For $\gamma\in\pi$, conjugates of $H_{\gamma,e}$ are of the form $H_{\delta,e}$ for some $\delta\in\pi$ by Lemma~\ref{lem:HgS operations}~\ref{enum:conjugation}. If $\gamma\neq \delta$, then $H_{\gamma,e}\cap H_{\delta,e}\in \F_1$ by Lemma~\ref{lem:HgS operations}~\ref{enum:intersection}. Hence the $\{H_{\gamma,e}\ |\ \gamma\in\pi\}$ are precisely the maximal elements in $\calD\setminus\F_1$ and condition $(M_{\F_1\subset \calD})$ holds.
		Moreover, given that $Z_\pi(\pi)$ is trivial, we have $N_{\pi\times\pi}(H_{\gamma,e})=H_{\gamma,e}$ by Lemma~\ref{lem:HgS operations}~\ref{enum:normalizer}.
	\end{proof}
\end{lem}

From now on and for the remainder of this note, we specialize to the following situation.

\begin{setup}\label{setup}
	Let $\pi$ be a torsionfree group admitting a malnormal collection of abelian subgroups $\calP=\{P_i\ |\ i\in I\}$ such that the centralizer $Z_\pi(b)$ is cyclic for any $b\in\pi$ that is not conjugate into any of the $P_i$.
\end{setup}
Note that in the situation of Setup~\ref{setup}, we have $N_\pi(Z_\pi(P_i))=Z_\pi(P_i)=P_i$ for every $P_i\in\calP$.
Our main examples of groups as in Setup~\ref{setup} are torsionfree relatively hyperbolic groups with finitely generated abelian peripheral subgroups, so-called \emph{toral} relatively hyperbolic groups.

The following lemma for the case when $\calP=\{\{e\}\}$ can be found in~\cite[Lemma 8.0.4]{Farber19} from where the first part of the proof is recalled. 
\begin{lem}\label{lem:centralizers}
	Let $\pi$ be a group as in Setup~\ref{setup}. Then for $b,c\in\pi\setminus\{e\}$, we have either $Z_\pi(b)=Z_\pi(c)$ or $Z_\pi(b)\cap Z_\pi(c)=\{e\}$.
	\begin{proof}
		Let $b,c\in\pi\setminus\{e\}$ be two elements. Suppose neither $b$ nor $c$ are conjugate into any of the $P_i$ and that $Z_\pi(b)\cap Z_\pi(c)$ is non-trivial. Let $Z_\pi(b)$, $Z_\pi(c)$ and $Z_\pi(b)\cap Z_\pi(c)$ be generated by $x$, $y$ and $z$, respectively. Then $x^n=z=y^m$ for some $n,m\in\IZ$. Observe that $z$ is not conjugate into any of the $P_i$. Thus its centralizer $Z_\pi(z)$ is infinite cyclic and contains both $x$ and $y$. Therefore, $x$ and $y$ commute and it follows that $Z_\pi(b)=Z_\pi(c)$.
		
		Suppose $b\in\pi\setminus\{e\}$ and $c\in gP_ig^{-1}$ for some $g\in\pi$, $P_i\in\calP$. Note that $Z_\pi(c)=gP_ig^{-1}$. If $Z_\pi(b)\cap gP_ig^{-1}$ is non-trivial, then $b\in gP_ig^{-1}$ by malnormality of $\calP$ and hence $Z_\pi(b)=Z_\pi(c)$. 
	\end{proof}
\end{lem}

\begin{lem}\label{lem:condition M}
	Let $\pi$ be a group as in Setup~\ref{setup}. Then we have the following:
	\begin{enumerate}[label=\enum]
		\item\label{enum:TRIV F1} Condition $(M_{\TR\subset \F_1})$ holds for the group $\pi\times\pi$. Moreover, for $\gamma\in\pi$ and $b\in \pi\setminus\{e\}$ there is an isomorphism $N_{\pi\times\pi}H_{\gamma,b}\cong Z_\pi(b)\times Z_\pi(b)$;
		\item\label{enum:TRIV F1restricted} Conditions $(M_{\TR\subset\F_1|_{H_{e,e}}})$ and $(NM_{\TR\subset\F_1|_{H_{e,e}}})$ hold for the group $H_{e,e}$.
	\end{enumerate}
	\begin{proof}
		\ref{enum:TRIV F1} For $\gamma\in\pi$ and $b\in\pi\setminus\{e\}$, conjugates of $H_{\gamma,b}$ are of the form $H_{\delta,c}$ for some $\delta\in\pi$, $c\in\pi\setminus\{e\}$ by Lemma~\ref{lem:HgS operations}~\ref{enum:conjugation}. We have either $H_{\gamma,b}=H_{\delta,c}$ or $H_{\gamma,b}\cap H_{\delta,c}=\{(e,e)\}$ by Lemma~\ref{lem:HgS operations}~\ref{enum:intersection} and Lemma~\ref{lem:centralizers}. Hence the $\{H_{\gamma,b}\ |\ \gamma\in\pi,b\in\pi\setminus\{e\}\}$ are precisely the maximal elements in $\F_1\setminus\TR$ and condition $(M_{\TR\subset \F_1})$ holds. 		
		Moreover, for $b\in\pi$ that is not conjugate into any of the $P_i$, observe that $N_\pi(Z_\pi(b))$ is torsionfree virtually cyclic and hence infinite cyclic. It follows that $N_\pi(Z_\pi(b))=Z_\pi(b)\cong\IZ$. If $b\in gP_ig^{-1}$ for some $g\in\pi$ and $P_i\in\calP$, we have $N_\pi(Z_\pi(b))=gP_ig^{-1}$ which is abelian and coincides with $Z_\pi(b)$.
		Thus, for any $b\in\pi\setminus\{e\}$ we have
		\[
			N_{\pi\times\pi}H_{\gamma,b}=\{(\gamma kh\gamma^{-1},h) \ |\ h,k\in Z_\pi(b)\} \cong Z_\pi(b)\times Z_\pi(b)
		\]
		by Lemma~\ref{lem:HgS operations}~\ref{enum:normalizer}.
				
		\ref{enum:TRIV F1restricted} Under the isomorphism $H_{e,e}\cong \pi$, the family $\F_1|_{H_{e,e}}$ is identified with the family $\F\spann{\{Z_\pi(b)\ |\ b\in\pi\setminus\{e\}\}}$. The claim follows as before by Lemma~\ref{lem:centralizers} and the observation that $Z_\pi(b)$ is self-normalizing for any $b\in\pi\setminus\{e\}$.
	\end{proof}
\end{lem}

\section{Maximality of topological complexity}

The following is the main technical result of this note and will immediately imply Theorem~\ref{thm:main thm} from the introduction.

\begin{thm}\label{thm:cd=2n}
	Let $\pi$ be a torsionfree group with $\cd(\pi)\ge 2$. Suppose that $\pi$ admits a malnormal collection of abelian subgroups $\calP=\{P_i\ |\ i\in I\}$ satisfying $\cd(P_i\times P_i)<\cd(\pi\times\pi)$ such that the centralizer $Z_\pi(b)$ is cyclic for any $b\in\pi$ that is not conjugate into any of the $P_i$. 	
	Then $\cd_{\TR\subset\calD}(\pi\times\pi)=\cd(\pi\times\pi)$.
	\begin{proof}
		We denote $\cd(\pi\times\pi)$ by $n$ and may assume that it is finite.
		Consider the families $\TR\subset \F_1\subset \calD$ of subgroups of $\pi\times\pi$ as defined in~\eqref{eqn:D F1}. 
		
		First, condition $(M_{\TR\subset \F_1})$ holds by Lemma~\ref{lem:condition M}~\ref{enum:TRIV F1} and hence Corollary~\ref{cor:LW}~\ref{enum:TRIV} yields a $(\pi\times\pi)$-pushout
		\begin{equation}\label{eqn:pushout1}\begin{tikzcd}
			\coprod_{H_{\gamma,b}\in\calM} (\pi\times\pi)\times_{N_{\pi\times\pi}H_{\gamma,b}}E(N_{\pi\times\pi}H_{\gamma,b}) \ar{r}\ar{d} & E(\pi\times\pi) \ar{d} \\
			\coprod_{H_{\gamma,b}\in\calM} (\pi\times\pi)\times_{N_{\pi\times\pi}H_{\gamma,b}}E(N_{\pi\times\pi}H_{\gamma,b}/H_{\gamma,b}) \ar{r} & E_{\F_1}(\pi\times\pi) \,,
		\end{tikzcd}\end{equation}
		where $\calM$ is a complete set of representatives of conjugacy classes of maximal elements in $\F_1\setminus\TR$. Moreover, in Lemma~\ref{lem:condition M}~\ref{enum:TRIV F1} we identified $N_{\pi\times\pi}H_{\gamma,b}\cong Z_\pi(b)\times Z_\pi(b)$ which is isomorphic to $\IZ\times\IZ$ or $P_i\times P_i$ for some $P_i\in\calP$ and hence has cohomological dimension strictly less than $n$. Thus, for any $\mathcal{O}_\calD(\pi\times\pi)$-module $M$, we have
		\[
			H^{n}_{\pi\times\pi}((\pi\times\pi)\times_{N_{\pi\times\pi}H_{\gamma,b}}E(N_{\pi\times\pi}H_{\gamma,b});M)=0 \,.
		\]
		Applying the Mayer--Vietoris sequence for $H^*_{\pi\times\pi}(-;M)$ to the pushout~\eqref{eqn:pushout1} yields that the map
		\[
			H^{n}_{\pi\times\pi}(E_{\F_1}(\pi\times\pi);M)\to H^{n}_{\pi\times\pi}(E(\pi\times\pi);M)
		\]
		is surjective.
		
		Second, conditions $(M_{\F_1\subset \calD})$ and $(NM_{\F_1\subset \calD})$ hold by Lemma~\ref{lem:F1 D} and hence Corollary~\ref{cor:LW}~\ref{enum:NM} yields a $(\pi\times\pi)$-pushout
		\begin{equation}\label{eqn:pushout2}\begin{tikzcd}
			(\pi\times\pi)\times_{H_{e,e}}E_{\F_1|_{H_{e,e}}}(H_{e,e})\ar{r}\ar{d} & E_{\F_1}(\pi\times\pi)\ar{d} \\
			(\pi\times\pi)/H_{e,e}\ar{r} & E_\calD(\pi\times\pi) \,.
		\end{tikzcd}\end{equation}
		Applying the Mayer--Vietoris sequence for $H^*_{\pi\times\pi}(-;M)$ to the pushout~\eqref{eqn:pushout2} yields that the map
		\[
			H^{n}_{\pi\times\pi}(E_\calD(\pi\times\pi);M)\to H^{n}_{\pi\times\pi}(E_{\F_1}(\pi\times\pi);M)
		\]
		is surjective provided that
		\begin{equation}\label{eqn:cohomology zero}
			H^{n}_{\pi\times\pi}((\pi\times\pi)\times_{H_{e,e}}E_{\F_1|_{H_{e,e}}}(H_{e,e});M) =0 \,.
		\end{equation}
		
		The latter is true by another application of Corollary~\ref{cor:LW}~\ref{enum:NM} using that conditions $(M_{\TR\subset\F_1|_{H_{e,e}}})$ and $(NM_{\TR\subset\F_1|_{H_{e,e}}})$ hold for the group $H_{e,e}$ by Lemma~\ref{lem:condition M}~\ref{enum:TRIV F1restricted}. It yields an $H_{e,e}$-pushout
		\begin{equation}\label{eqn:pushout3}\begin{tikzcd}
			\coprod_{H_{e,b}\in\calM'}H_{e,e}\times_{H_{e,b}}E(H_{e,b})\ar{r}\ar{d} & E(H_{e,e})\ar{d} \\
			\coprod_{H_{e,b}\in\calM'}H_{e,e}/H_{e,b}\ar{r} & E_{\F_1|_{H_{e,e}}}(H_{e,e}) \,,
		\end{tikzcd}\end{equation}
		where $\calM'$ is a complete set of representatives of conjugacy classes of maximal elements in $\F_1|_{H_{e,e}}\setminus\TR$.
		The Mayer--Vietoris sequence for $H^*_{H_{e,e}}(-;M)$ applied to the pushout~\eqref{eqn:pushout3} shows that~\eqref{eqn:cohomology zero} indeed holds, using that $\cd(H_{e,e})<n$ and $\cd(H_{e,b})<n-1$ for $b\in\pi\setminus\{e\}$. 
	
		Together, the map
		\[
			H^{n}_{\pi\times\pi}(E_\calD(\pi\times\pi);M)\to H^{n}_{\pi\times\pi}(E(\pi\times\pi);M)
		\]
		is surjective for any $\mathcal{O}_\calD(\pi\times\pi)$-module $M$. Finally, the coefficients $M$ can be chosen such that $H^{n}_{\pi\times\pi}(E(\pi\times\pi);M)$ is non-trivial.
		This concludes the proof.
	\end{proof}
\end{thm}

\paragraph{Proof of Theorem~\ref{thm:main thm}}
It follows from Theorem~\ref{thm:cd=2n} that the inequalities
\[
	\cd_{\TR\subset\calD}(\pi\times\pi)\le \TC(\pi)\le \cd(\pi\times\pi)
\]
are in fact equalities.\qed

\bibliographystyle{alpha}
\bibliography{bib}

\begin{thebibliography}{FGLO19}

\bibitem[BCE]{BCE20}
Zbigniew B\l{}aszczyk, Jos\'{e} Carrasquel, and Arturo Espinosa.
\newblock On the sectional category of subgroup inclusions and {A}damson
  cohomology theory.
\newblock Preprint, arXiv:2012.11912, 2020.

\bibitem[CLM]{CLM20}
Pietro Capovilla, Clara L\"{o}h, and Marco Moraschini.
\newblock Amenable category and complexity.
\newblock Preprint, arXiv:2012.00612, 2020.

\bibitem[Dra19]{Dranishnikov19}
Alexander Dranishnikov.
\newblock On dimension of product of groups.
\newblock {\em Algebra Discrete Math.}, 28(2):203--212, 2019.

\bibitem[Dra20]{Dranishnikov20}
Alexander Dranishnikov.
\newblock On topological complexity of hyperbolic groups.
\newblock {\em Proc. Amer. Math. Soc.}, 148(10):4547--4556, 2020.

\bibitem[Far03]{Farber03}
Michael Farber.
\newblock Topological complexity of motion planning.
\newblock {\em Discrete Comput. Geom.}, 29(2):211--221, 2003.

\bibitem[FGLO19]{Farber19}
Michael Farber, Mark Grant, Gregory Lupton, and John Oprea.
\newblock Bredon cohomology and robot motion planning.
\newblock {\em Algebr. Geom. Topol.}, 19(4):2023--2059, 2019.

\bibitem[Flu]{Fluch10}
Martin Fluch.
\newblock On {B}redon (co-)homological dimension of groups.
\newblock PhD thesis, arXiv:1009.4633, 2010.

\bibitem[FM20]{Farber-Mescher20}
Michael Farber and Stephan Mescher.
\newblock On the topological complexity of aspherical spaces.
\newblock {\em J. Topol. Anal.}, 12(2):293--319, 2020.

\bibitem[L{\"u}c05]{Lueck05survey}
Wolfgang L{\"u}ck.
\newblock Survey on classifying spaces for families of subgroups.
\newblock In {\em Infinite groups: geometric, combinatorial and dynamical
  aspects}, volume 248 of {\em Progr. Math.}, pages 269--322. Birkh\"{a}user,
  Basel, 2005.

\bibitem[LW12]{Lueck-Weiermann12}
Wolfgang L\"{u}ck and Michael Weiermann.
\newblock On the classifying space of the family of virtually cyclic subgroups.
\newblock {\em Pure Appl. Math. Q.}, 8(2):497--555, 2012.

\end{thebibliography}

\end{document}